\documentclass[11pt, reqno]{amsart}
\usepackage{mathrsfs}
\usepackage{graphics}
\usepackage{latexsym,amssymb,amsmath,amscd,amscd,amsthm,amsxtra}
\usepackage[dvips]{graphicx}
\usepackage[dvips]{graphics,color}
\usepackage{amsfonts}
   \def\lineno.sty{\texttt{\itshape lineno.sty}}
\def\Box{\vcenter{\vbox{\hrule\hbox{\vrule
     \vbox to 8.8pt{\hbox to 10pt{}\vfill}\vrule}\hrule}}}

\newcommand{\F}{ {{\mathbb F}} }

\def\Im{\operatorname{Im}}

\newtheorem{thm}{Theorem}[section]
\newtheorem{lem}[thm]{Lemma}
\newtheorem{cor}[thm]{Corollary}

\newtheorem{prop}[thm]{Proposition}

\newtheorem{remark}{Remark}

\numberwithin{equation}{section}

\begin{document}
 \thispagestyle{empty}
\title{Linearized Wenger graphs}

\author[X. Cao, M. Lu, D. Wan, L.-P. Wang, Q. Wang]{ Xiwang Cao, Mei Lu, Daqing Wan, Li-Ping Wang, Qiang Wang}

\address{
Xiwang Cao is with the School of Mathematical Sciences, Nanjing University of
Aeronautics and Astronautics, Nanjing 210016, China, email: {\tt
xwcao@nuaa.edu.cn}}
\address{Mei Lu is with Department of Mathematical Sciences, Tsinghua
University, Beijing 100084, China, email: {\tt
mlu@math.tsinghua.edu.cn}}
\address{Daqing Wan is with Department of Mathematics, University of California, Irvine, CA 92697-3875, USA, email: {\tt dwan@math.uci.edu}}
\address{Li-Ping Wang is with Institute of Information Engineering, Chinese Academy of Sciences, Beijing 100093, Beijing, China, email: {\tt wangliping@iie.ac.cn}}
\address{Qiang Wang is with School of Mathematics and Statistics, Carleton University, 1125 Colonel By Drive, Ottawa, Ontario K1S 5B6, Canada. email: {\tt
wang@math.carleton.ca}}
\thanks{The work of this paper was supported by National Natural Science Foundation under grant number 11371011 and No. 61170289 and China Scholarship Council.}

\begin{abstract}
Motivated by recent extensive studies on Wenger graphs, we introduce a new infinite class of bipartite graphs of the similar type, called linearized Wenger graphs. The spectrum, diameter and girth of these linearized Wenger graphs are determined.
\end{abstract}

\keywords{Cayley graph, graph spectrum, expander, algebraic graph theory, diameter, girth}

\maketitle

\section{Introduction }

Let $\mathbb{F}_q$ be a finite field of order $q$ such that $p$ is prime and $q=p^e$ a prime power. All graph theory notions can be found in Bollob\'{a}s \cite{boll}.
Recently, a class of bipartite graphs called {\it Wenger graphs} which are defined over $\mathbb{F}_q$ has attracted a lot of attention  because of their nice graphical properties  \cite{cll, lu, lu2, MM,SHS,Viglione0,Viglione,Wenger}.  For example, the number of edges of these graphs meets the lower bound of Tur\'{a}n number of the cycle with length $4, 6, 10$  \cite{Wenger}. The original definition was introduced by Wenger \cite{Wenger} for $p$-regular bipartite graphs and then was extended by Lazbnik and Ustimenko \cite{lu} for arbitrary prime power $q$.  An equivalent representation of these graphs appeared later in Lazebnik and Viglione \cite{LV}  and then a more general class of graphs was defined in \cite{Viglione0}, on which we concentrate in this paper.

Let $m\geq 1$ be a positive integer and $g_k(x,y)\in \mathbb{F}_q[x,y]$ for $ 2\leq k \leq m+1$. Let $\mathfrak{P}=\mathbb{F}_q^{m+1}$ and  $\mathfrak{L}=\mathbb{F}_q^{m+1}$ be two copies of the $(m+1)$-dimensional vector space over  $\F_q$, which are called the point set and the line set respectively.
Let $\mathfrak{G} = G_q(g_2, \cdots, g_{m+1}) =(V,E)$ be the graph with vertex set $V=\mathfrak{P}\cup \mathfrak{L}$ and the edge set $E$ is defined as follow: there is an edge from a point $P=(p_1,p_2,\cdots,p_{m+1})\in \mathfrak{P}$  to a line $L=[l_1,l_2,\cdots,l_{m+1}]\in \mathfrak{L}$, denoted by $P\sim L$ (we force $\mathfrak{G}$ to be a undirected graph by removing the arrows),  if the following $m$ equalities hold:
\begin{eqnarray}\label{defnGraph}
  l_2+p_2 &=& g_2(p_1,l_1)  \nonumber  \\
 l_3+p_3 &=& g_3(p_1,l_1) \nonumber  \\
 \vdots&\vdots& \vdots\\
 l_{m+1}+p_{m+1} &=& g_{m+1}(p_1,l_1).  \nonumber
\end{eqnarray}

 If $g_k(x, y), k=2, \cdots, m+1$, are all monomials, the graph is called a {\it monomial graph}; see  \cite{DLW}.
 If $g_k(x,y)=x^{k-1}y, k=2,\cdots,m+1$, then the graph is just the original Wenger graph in \cite{cll}, also denoted by $W_m(q)$.
It was shown in \cite{lu} that the automorphism group of $W_m(q)$ acts transitively on each of $\mathfrak{P}$ and $\mathfrak{L}$, and on the set of edges of $W_m(q)$. In other words, the graphs $W_m(q)$ are point-, line-, and edge-transitive. It is also shown that, see \cite{lu2}, $W_1(q)$ is vertex-transitive for all $q$, and that $W_2(q)$ is vertex-transitive for even $q$. For all $m\geq 3$ and $q\geq 3$, and for $m=2$ and all odd $q$, the graphs $W_m(q)$ are not vertex-transitive. Another result of \cite{lu2} is that $W_m(q)$ is connected when $1\leq m\leq q-1$, and disconnected when $m\geq q$, in which case it has $q^{m-q+1}$ components, each isomorphic to $W_{q-1}(q)$. In \cite{Viglione}, Viglione proved that  the diameter of $W_m(q)$ is $2m+2$ when $1\leq m\leq q-1$.  In \cite{cll}, Cioab$\rm{\breve{a}}$, Lazebnik and Li determined the spectrum of $W_m(q)$.

\medskip
In this paper we focus on the basic properties of some extensions of Wenger graphs defined as in Equation~(\ref{defnGraph}). In Section 2  we first study the spectrum of a general class of graphs such that polynomials $g_k(x,y)\in \mathbb{F}_q[x,y]$ are defined by $g_k(x,y)=f_k(x)y$, and the mapping $\vartheta: \mathbb{F}_q\rightarrow \mathbb{F}_q^{m+1}; u\mapsto (1,f_2(u),\cdots,f_{m+1}(u))$ is injective. The eigenvalues of such a graph are determined, however, their multiplicities are reduced to counting
certain polynomials with a given number of roots over finite fields. The latter problem is an interesting number theoretical
problem, which is expected to be difficult in general. A complete solution in interesting special cases is already significant.
In particular, we introduce a new class of bipartite graphs called  linearized Wenger graphs. These graphs are denoted by $L_m(q)$, which are defined by Equation~(\ref{defnGraph}) together with  $g_k(x,y)=x^{p^{k-2}}y, k=2,\cdots,m+1$. Using results on linearized polynomials over finite fields, we are able to explicitly determine the spectrum of such graphs when $m\geq e$ in Section 3. Finally we obtain the diameter and girth of linearized Wenger graphs in Section 4 and Section 5, respectively. As a consequence, when $m=e$,  this provides a new class of infinitely many connected $p^e$-regular expander graphs of $q^{2m+2}$ vertices with optimal diameter $2(m+1)$ when either the prime $p$ or the exponent $e$ goes to infinity.

\section{The spectrum of general Wenger graphs}

In this section we study the basic properties of the class of graphs $\mathfrak{G}$ defined by $g_k(x,y)=f_k(x)y$, where $g_k(x, y)$ is a product of a polynomial in terms of $x$ and the linear polynomial $y$, for $2\leq k \leq m+1$.

{\prop The graph $\mathfrak{G}=G_q(f_2(x)y,\ldots,f_{m+1}(x)y)$ is $q$-regular.}
\begin{proof}
Given a point $P$ and a line $L$ in $V$, by definition, $P=(p_1,p_2,\cdots,p_{m+1})$ is adjacent to $L=[l_1,l_2,\cdots,l_{m+1}]$ if and only if the following $m$ equalities hold:
\begin{equation}\label{defnGraph1}
\left\{\begin{array}{ccc}
  l_2+p_2 &=& f_2(p_1)l_1    \\
 l_3+p_3 &=& f_3(p_1)l_1   \\
 \vdots&\vdots& \vdots\\
 l_{m+1}+p_{m+1} &=& f_{m+1}(p_1)l_1.
       \end{array}
\right.
\end{equation}
When the point $P$ is prescribed, (\ref{defnGraph1}) implies that one can uniquely solve $l_k$ ($k\geq 2$) from $l_1$, and thus (\ref{defnGraph1}) has $q$ solutions.
Similarly, when the point $L$ is prescribed, (\ref{defnGraph1}) implies that one can uniquely solve $p_k$ ($k\geq 2$)  from $p_1$, and thus (\ref{defnGraph1}) has $q$ solutions.
\end{proof}

Since $\mathfrak{G}$ is a bipartite graph, its adjacency matrix  is of the form:
\begin{equation*}
  A=\left(\begin{array}{cc}
            0 & N \\
            N^T & 0
          \end{array}
  \right)
\end{equation*}
with a matrix $N$ and
\begin{equation}\label{f-4'}
  A^2=\left(\begin{array}{cc}
           NN^T & 0\\
          0& N^TN
          \end{array}
  \right).
\end{equation}

In order to consider the properties of $\mathfrak{G}$, we define a graph $H$ as follows:  the vertex set is $\mathbb{F}_q^{m+1}$ containing all lines in $\mathfrak{G}$,  any two lines $L=[l_1,l_2,\cdots,l_{m+1}]$ and $L'=[l_1',l_2',\cdots,l_{m+1}']$ are adjacent if and only if they share a common neighbor point  $P=(p_1,p_2,\cdots, p_{m+1})$ in the graph $\mathfrak{G}$ defined above.

Moreover, one can check that the graph $H$ is a Cayley graph  with the generating set \begin{equation*}
  S=\{(t,tf_2(u),\cdots,tf_{m+1}(u))|\, t\in \mathbb{F}_q^*,u\in \mathbb{F}_q\}.
\end{equation*}
Indeed, $L\sim L'$ if and only if
$l_k -l_k' = f_k(p_1, l_1) - f_k(p_1, l_1') = f_k(p_1) (l_1 -l_1')$ for $2\leq k\leq m+1$.

Furthermore,  if $B$ is the adjacency matrix of $H$  then
\begin{equation}\label{f-4}
 NN^T=B+qI,
\end{equation}
where $I$ is the identity matrix.  Let us denote all eigenvalues of $H$ by $\lambda_1(B)$, $\ldots$, $\lambda_{q^{m+1}}(B)$ .  Since $N^TN$ and $NN^T$ have the same eigenvalues,  one can check that the eigenvalues of $\mathfrak{G}$ are $\pm\sqrt{\lambda_i(B)+q}, i=1,2,\cdots,q^{m+1}$.

Now let us assume the mapping $\vartheta: \mathbb{F}_q\rightarrow \mathbb{F}_q^{m+1}; u\mapsto (1,f_2(u),\cdots,f_{m+1}(u))$ is injective.  Then we know that $|S|=q(q-1)$.
Our first result is the following

{\thm \label{thm-1} Let $\mathfrak{G}$ be defined in (\ref{defnGraph})
 with the assumptions that  $g_k(x,y)=f_k(x)y$ for $k=2,\cdots,m+1$
and the mapping $\vartheta: \mathbb{F}_q\rightarrow \mathbb{F}_q^{m+1}$ defined by
$u\mapsto (1,f_2(u),\cdots,f_{m+1}(u))$ is injective.
For all prime power $q$ and positive integer $m$, the eigenvalues of $\mathfrak{G}$, counted with multiplicities,  are
$$\pm \sqrt{qN_{F_w}}, w=(w_1,w_2,\cdots,w_{m+1})\in \mathbb{F}_q^{m+1},$$
where $F_w(u)=w_1+w_2f_2(u)+\cdots+w_{m+1}f_{m+1}(u)$ and $N_{F_w}=|\{u\in \mathbb{F}_q: F_w(u)=0\}|$. For $0\leq i\leq q$, the multiplicity of $\pm \sqrt{q{i}}$ is
$$n_i=|\{w\in \mathbb{F}_q^{m+1}:N_{F_w}=i\}|.$$
Moreover, the number of connected components of $\mathfrak{G}$ is
\begin{equation*}
q^{m+1-{\rm rank}_{\mathbb{F}_q}(1,f_2,\cdots,f_{m+1})}.
\end{equation*}
Therefore $\mathfrak{G}$ is connected if and only if
$1,f_2,\cdots,f_{m+1}$ are $\mathbb{F}_q$-linearly independent.
}
\begin{proof}
Let $\zeta_p$ be a primitive $p$-th root of unity, and for every $w :=(w_1,w_2,\cdots,w_{m+1})\in \mathbb{F}_q^{m+1}$, we define a character $\psi_w: \mathbb{F}_q^{m+1}\rightarrow \mathbb{C}^*$ by
\begin{equation*}
 \psi_w:  u=(u_1,u_2,\cdots,u_{m+1})\mapsto \zeta_p^{{\rm tr}(w_1u_1+w_2u_2+\cdots+w_{m+1}u_{m+1})},
\end{equation*}
where ${\rm tr}$ is the absolute trace map. As described in \cite{babai,lovasz},
the eigenvalues of the Cayley graph $H$ are
\begin{equation}\label{f-6}
   \psi_w(S):=\sum_{t\in \mathbb{F}_q^*,u\in \mathbb{F}_q}\zeta_p^{{\rm tr}(t(w_1+w_2f_2(u)+\cdots+w_{m+1}f_{m+1}(u)))}, w\in \mathbb{F}_q^{m+1}.
\end{equation}
Denote by $F_w(u)$ the function $w_1+w_2f_2(u)+\cdots+w_{m+1}f_{m+1}(u)$ and
$N_{F_w}=|\{u\in \mathbb{F}_q: F_w(u)=0\}|$. Then it follows that
\begin{eqnarray*}
  \psi_w(S)&= & \sum_{t\in \mathbb{F}^*_q,u\in \mathbb{F}_q}\zeta_p^{{\rm tr}(tF_w(u))}\\
           &= & \sum_{t\in \mathbb{F}^*_q, F_w(u)=0}\zeta_p^{{\rm tr}(tF_w(u))}+
           \sum_{t\in \mathbb{F}^*_q, F_w(u)\neq 0}\zeta_p^{{\rm tr}(tF_w(u))} \\
           &=&(q-1) N_{F_w} + (-1) (q - N_{F_w}) \\
&=& q\left(N_{F_w}-1\right).
\end{eqnarray*}
Thus this derives  that the eigenvalues of $\mathfrak{G}$ are
\begin{equation}\label{f-7}
 \pm \sqrt{qN_{F_w}}, w\in \mathbb{F}_q^{m+1},
\end{equation}
where $N_{F_w}=|\{u\in \mathbb{F}_q: F_w(u)=0\}|$.  For example, when $w=(0, \ldots, 0)$ we have $N_{F_0} = q$ which implies that $\mathfrak{G}$ has $\pm q$ as its eigenvalues. Moreover, for any $w\neq 0$, it is easy to see that $N_{F_w} \leq \deg(F_w) \leq \max\{\deg(f_2), \ldots, \deg(f_{m+1})\}$.

The number of connected components of $\mathfrak{G}$ is
\begin{equation}\label{f-8}
  |\{w: F_w(x)\equiv 0\mbox{ for all $x\in \mathbb{F}_q$}\}|=q^{m+1-{\rm rank}_{\mathbb{F}_q}(1,f_2,\cdots,f_{m+1})}.
\end{equation}
Therefore $\mathfrak{G}$ is connected if and only if
$1,f_2,\cdots,f_{m+1}$ are $\mathbb{F}_q$-linearly independent.
\end{proof}

\begin{remark} The computation of the multiplicities $n_i$'s is obviously an interesting number theoretical problem.
One cannot expect a simple closed formula for $n_i$'s in general. Among the most interesting case is when the
$f_k(x)$'s are given by monomials in $x$. When the $f_k$'s are consecutive monomials (the original Wenger
graph), there is indeed a simple formula for $n_i$'s. When the $f_k$'s  are not consecutive monomials,
the problem is more difficult. The linearized Wenger graph considered in next section deals with the first
non-trivial example of non-consecutive monomials.
\end{remark}

\section{The spectrum of  linearized Wenger graphs}

Let $q=p^e$ and $m$ be a positive integer as before. We focus on the linearized Wenger graph $L_m(q)$ from now on where $f_k(x)=x^{p^{k-2}}$,  $k=2,\cdots,m+1$. The goal of this section is to explicitly compute the  spectum  of $L_m(q)$ by determining the explicit formula of $N_{F_w}$ and $n_i$
in Theorem \ref{thm-1}.  The computation involved in linearized Wenger graphs is more complicated  since the degrees of $f_k(x)=x^{p^{k-2}}, k=2,\ldots,m+1$  are high and not consecutive as in Wenger graphs.

We first give a basic lemma which will be used in the rest of the paper. It is an old result with
the first derivation of the formula due to Landsberg [9, p.455]; see also Lemma 2.1 in \cite{ling}.

{\lem \label{lem-3.1}
The number of $l \times n$ matrices  over $\mathbb{F}_q$ with rank $k$ is
 $ \frac{\prod_{i=0}^{k-1}(q^l-q^i)(q^n-q^i)}{\prod_{i=0}^{k-1}(q^k-q^i)}$. }

\begin{proof}
For a fixed $k$-dimensional subspace $W\in \mathbb{F}^{l}_{q}$, the number of $l\times n$ matrices
with $W$ as the column space is equal to the number of $k\times n$ matrices of rank $k$. Such
a matrix is given by the $k$ linearly independent row vectors of length $n$. The number
of those is $\prod_{i=0}^{k-1}(q^n-q^i)$. The number of $k$-dimensional subspaces of $\mathbb{F}^{l}_{q}$
is $\frac{\prod_{i=0}^{i}(q^l-q^i)}{\prod_{i=0}^{i}(q^k-q^i)}$ and the product is the number of rank $k$ matrices.
\end{proof}

When $m= e$, the functions $1,x,\cdots,x^{p^{m-1}}$ are $\mathbb{F}_q$-linearly independent and so $L_m(q)$ is connected. For every $w=(w_1,w_2,\cdots,w_{m+1})\in \mathbb{F}_q^{m+1}$, define $F_w(x)=w_1+w_2x +w_3x^{p} +\cdots+w_{m+1}x^{p^{m-1}}$. By Theorem \ref{thm-1}, the eigenvalues of the linearized Wenger graph $L_m(q)$, counting multiplicities,  are
\begin{equation*}
  \pm \sqrt{qN_{F_w}}, w\in \mathbb{F}_q^{m+1},
\end{equation*}
where $N_{F_w}=|\{u\in \mathbb{F}_q: F_w(u)=0\}| =   |\{u\in \mathbb{F}_q: \bar{F}_w(u)=-w_1\}|$,
where  $\bar{F}_w(x)= w_2x+\cdots+w_{m+1}x^{p^{m-1}}$ is an $\mathbb{F}_p$-linearized polynomial. If $-w_1 \not\in \Im(\bar{F}_w)$, then $N_{F_w} =0$. Otherwise, this also implies that
\begin{equation*}
  N_{F_w}=p^{\dim_{\mathbb{F}_p}(\ker(\bar{F}_w))}.
\end{equation*}

Choosing a fixed basis of $\mathbb{F}_q/\mathbb{F}_p$ as $\alpha_1,\cdots,\alpha_e$, we know that every $p$-linear polynomial $\bar{F}_w(x)$ can be written as
\begin{equation}\label{f-10}
 \bar{F}_w(x)={\rm tr}(\beta_1x)\alpha_1+{\rm tr}(\beta_2x)\alpha_2+\cdots+{\rm tr}(\beta_ex)\alpha_e,
\end{equation}
where $\beta_1,\cdots,\beta_e$ are elements in $\mathbb{F}_q$ uniquely determined by
$w_2, \ldots, w_{m+1}$.  By Theorem~2.2 in \cite{ling}, we have $\dim_{\mathbb{F}_p}(\ker(\bar{F}_w))=i$ if and only if ${\rm rank}_{\mathbb{F}_p}(\beta_1,\cdots,\beta_e) = e-i$. For $0\leq i \leq e$, there are exactly
\[
\frac{\prod_{j=0}^{e-i-1} (p^e - p^j)^2}{ \prod_{j=0}^{e-i-1} (p^{e-i} - p^j)}
\]
different $w_2, \ldots, w_{m+1}$ such that $\dim_{\mathbb{F}_p}(\ker(\bar{F}_w)) =i$ by Lemma \ref{lem-3.1}.  There are $p^{e-i}$ choices for $-w_1$ in the image set of $\bar{F}_w$, therefore the multiplicity of the eigenvalue $\pm \sqrt{q p^i}$ is
\begin{equation}\label{f-9}
  n_{p^i}= p^{e-i} \frac{\prod_{j=0}^{e-i-1} (p^e - p^j)^2}{ \prod_{j=0}^{e-i-1} (p^{e-i} - p^j)}.
\end{equation}

Now, counting each $-w_1$ not in the image set of $\bar{F}_w$ such that $\dim_{\mathbb{F}_p}(\ker(\bar{F}_w)) =i$ for $ 1\leq i\leq e$,  the multiplicity of the eigenvalue $0$ is
\begin{equation}\label{eigenvalue0}
  n_0= \sum_{i=1}^{e} (p^e-p^{e-i}) \frac{\prod_{j=0}^{e-i-1} (p^e - p^j)^2}{ \prod_{j=0}^{e-i-1} (p^{e-i} - p^j)}.
\end{equation}






When $m>e$, one checks that ${\rm rank}_{\mathbb{F}_q}(1,x,x^p, \cdots, x^{p^{m-1}})=e+1$ and thus we obtain the following result:

{\thm \label{thm-3}Let $m\geq e$. The linearized Wenger graph $L_m(q)$ has $q^{m-e}$ components. The distinct eigenvalues are
\begin{equation*}
  ~ 0, ~\pm \sqrt{qp^i}, 0\leq i\leq e.
\end{equation*}
For $0\leq i\leq e$, the multiplicity of the eigenvalue $\pm \sqrt{qp^i}$ is  $ q^{m-e} n_{p^i}$ where $n_{p^i}$ is given by (\ref{f-9}). The multiplicity of the eigenvalue $0$ is
$q^{m-e} n_0$ where $n_0$ is given by  (\ref{eigenvalue0}). }

When $m=e$, these linearized Wenger graphs are connect $q$-regular $(q, \epsilon)$-expander graphs with edge expansion $\epsilon > \frac{q - \sqrt{qp^{e-1}}}{2} = \frac{q^{1/2} p^{(e-1)/2}(p^{1/2} -1)}{2}$. As to expander graphs, we refer to \cite{HLW,jukna} for more details.

When $m<e$, the linearized Wenger graph $L_m(q)$ is connected, however, we do not know a closed formula for the
multiplicities of the eigenvalues $\pm \sqrt{qp^i}$.  We leave this as an open problem.

\section{The diameter of linearized Wenger graphs}

 Recall that a sequence of vertices $v_1,\cdots,v_s$ in a simple graph  $\mathfrak{G}=(V,E)$ defines a {\it path} of length $s-1$ if $(v_i,v_{i+1})\in E$ for every $i, 1\leq i\leq s-1$. The {\em distance} between $v_i$ and $v_j$ is the number of edges in a shortest path joining $v_i$ and $v_j$. The {\em diameter} of a graph $\mathfrak{G}$ is the maximum distance between any two vertices of $\mathfrak{G}$. In  \cite{Viglione} it is shown that the diameter of the Wenger graph $W_m(q)$ is $2m+2$ when $1\leq m\leq q-1$. In this section, we assume that $m\leq e$ so that the linearized Wenger graphs are connected.  We now explicitly determine  the diameter of the linearized Wenger graph $L_m(q)$.

{\thm \label{thm-5} If $m\leq e$, the diameter  of the linearized Wenger graph $L_m(q)$ is $2(m+1)$.}

Before proceeding to the proof of the above theorem, we give the following lemma.

{\lem \label{lem-4.1}If $x_1,\ldots, x_{m}$ in $\mathbb{F}_q$ are $\mathbb{F}_p$-linearly independent, then
$$\left |\begin{array}{cccc}
             1  & 1 & \ldots  & 1 \\
             x_1 & x_2 & \ldots & x_{m}\\
             {x}^{p}_1 & {x}^{p}_2 & \ldots & {x}^{p}_m \\
             \vdots & \vdots & \vdots & \vdots \\
             x^{p^{m-2}}_1 & x^{p^{m-2}}_2 & \ldots & x^{p^{m-2}}_m
           \end{array}\right |\neq 0.$$}


\begin{proof}
First it is easy to see that
$$\left |\begin{array}{cccc}
             1  & 1 & \ldots  & 1 \\
             x_1 & x_2 & \ldots & x_{m}\\
             {x}^{p}_1 & {x}^{p}_2 & \ldots & {x}^{p}_m \\
             \vdots & \vdots & \vdots & \vdots \\
             x^{p^{m-2}}_1 & x^{p^{m-2}}_2 & \ldots & x^{p^{m-2}}_m
           \end{array}\right |=\left |\begin{array}{cccc}
             1  & 1 & \ldots  & 1 \\
             0 & x_2-x_1 & \ldots & x_{m}-x_1\\
             0 & (x_2-x_1)^{p} & \ldots & (x_m-x_1)^{p} \\
             \vdots & \vdots & \vdots & \vdots \\
             0 & (x_2-x_1)^{p^{m-2}} & \ldots & (x_m-x_1)^{p^{m-2}}
           \end{array}\right |.$$

\noindent Since $x_1,\ldots,x_m$ are $\mathbb{F}_p$-linearly independent, $x_2-x_1,\ldots, x_{m}-x_1$ are $\mathbb{F}_p$-linearly independent. By induction,
 $\left |\begin{array}{ccc}
              x_2-x_1 & \ldots & x_{m}-x_1\\
              (x_2-x_1)^{p} & \ldots & (x_m-x_1)^{p} \\
              \vdots & \vdots & \vdots \\
              (x_2-x_1)^{p^{m-2}} & \ldots & (x_m-x_1)^{p^{m-2}}
           \end{array}\right |\neq 0$, the proof is complete.
\end{proof}

\begin{proof}[Proof of Theorem \ref{thm-5}]

First we consider the distance between any two vertices $L$ and $L'$ in $\mathfrak{L}$  of the linearized Wenger graph $L_m(q)$. If $L_1P_1\ldots P_sL_{s+1}$ is a path in $L_m(q)$ between  $L=L_1$ and $L'=L_{s+1}$, where $L_i=[l_1^{(i)},\cdots,l_{m+1}^{(i)}]$ and
$P_i=(p_1^{(i)},\cdots,p_{m+1}^{(i)})$, we have
\begin{equation*}
  l_k^{(i+1)}-l_k^{(i)}=(l_1^{(i+1)}-l_1^{(i)})(p_1^{(i)})^{p^{k-2}}, k=2,\cdots,m+1, i=1,\cdots,s.
\end{equation*}
 Therefore  there are elements $t_i =l_1^{(i+1)}-l_1^{(i)} $, $x_i = p_1^{(i)} \in \mathbb{F}_q$, $1\leq i\leq s$,  such that
\begin{equation}\label{f-11}
  (L_{s+1}-L_1)^T =t_1\left(\begin{array}{c}
             1 \\
             x_1 \\
             {x}^{p}_1 \\
             \vdots\\
             x^{p^{m-1}}_1
           \end{array}
  \right)+t_2\left(\begin{array}{c}
             1 \\
             x_2 \\
             {x}^{p}_2 \\
             \vdots\\
             {x}^{p^{m-1}}_2
           \end{array}
  \right)+\cdots+t_s\left(\begin{array}{c}
             1 \\
             x_s \\
             {x}^{p}_s \\
             \vdots\\
             {x}^{p^{m-1}}_s
           \end{array}
  \right ).
\end{equation}


Take $s=m+1$ and choose $x_1, \ldots, x_{m+1} \in \mathbb{F}_q$ such that  $x_2-x_1,\ldots,x_{m+1}-x_1$ are $\mathbb{F}_p$-linearly independent. Then by Lemma~\ref{lem-4.1},  the  coefficient matrix of  Eq. (\ref{f-11}) is nonsingular, and  thus Eq. (\ref{f-11}) has a unique solution for $t_1, t_2, \ldots, t_s$. Thus the distance of any two vertices in $\mathfrak{L}$  is at most  $2(m+1)$.

Similarly,  let us consider any two vertices $P$ and $P'$ in $\mathfrak{P}$ of $L_m(q)$. Let $P_1L_1\ldots L_sP_{s+1}$ is a path in $L_m(q)$ between  $P=P_1$ and $P'=P_{s+1}$, where $L_i=[l_1^{(i)},\cdots,l_{m+1}^{(i)}]$ and $P_i=(p_1^{(i)},\cdots,p_{m+1}^{(i)})$.  Then we have
\begin{equation*}
  p_k^{(i+1)}-p_k^{(i)}=l_1^{(i)}(p_1^{(i+1)}-p_1^{(i)})^{p^{k-2}}, k=2,\cdots,m+1, i=1,\cdots,s.
\end{equation*}
Similarly, if we take $s= m+1$ and choose $p_i \in \F_q$ such that $p_1^{(i+1)}-p_1^{(i)}$, $1 \leq i \leq m$ are $\F_p$-linearly independent,  then we can find unique solution for $l_1^{(1)}, \ldots, l_1^{(m)}$. Hence  the distance of any two vertices in $\mathfrak{P}$  is at most $2(m+1)$.

 Finally, we consider the distance between a vertex $P=(p_1,\ldots,p_{m+1})\in \mathfrak{P}$ and a vertex $L \in \mathfrak{L}$.  First we choose any line $L_1\in \mathfrak{L}$ such that it is adjacent to $P$.  From the earlier discussion,  there exists a path  from $L_1$ to $L$ with distince at most $2(m+1)$. We modify the earlier construction so that the path goes through the vertex $P$. Namely,  In Eq. (\ref{f-11}), we let
 $x_1=p_1$ and choose the rest of $x_i$'s so that $x_2-x_1,\ldots, x_{m+1}-x_{1}\in \mathbb{F}_q$ are $\mathbb{F}_p$-linearly independent. Then there is a unique solution $\{t_1,\ldots,t_{s}\}$ and so there is a path between $L_1$ and $L$ with length at most $2(m+1)$ passing through $P$. Therefore the distance of $P$ and $L$ is less than or equal to $2(m+1)$.  Hence the diameter of $L_m(q)$ is always at most $2(m+1)$.

On the other hand, we now show that the distance $2(m+1)$ can be reached. Indeed, choose two vertices $L_1$ and $L_{s+1}$ such that $L_{s+1}-L_1=[0,\ldots,0,1]$. We can show that the distance between them is at least $2(m+1)$.
Otherwise, suppose there is a path from $L_1$ to $L_{s+1}$ with distance $2s \leq 2m$.  Then Eq. (\ref{f-11}) has a solution with $1\leq s\leq m$. We show that this is impossible.

If  either $x_1,\ldots,x_{s}$ are $\mathbb{F}_p$-linearly independent and $s<m$, or $x_1,\ldots,x_{s}$ are $\mathbb{F}_p$-linearly dependent,  then the last $m$ rows of (\ref{f-11}) always can be reduced to
\begin{equation}\label{f-12}
 \left(\begin{array}{c}
             0 \\
             0 \\
             \vdots\\
             1
           \end{array}
  \right) =t'_1\left(\begin{array}{c}
             x'_1 \\
             (x'_1)^{p} \\
             \vdots\\
             (x'_1)^{p^{m-1}}
           \end{array}
  \right)+t'_2\left(\begin{array}{c}
             x'_2 \\
             (x'_2)^{p} \\
             \vdots\\
             (x'_2)^{p^{m-1}}
           \end{array}
  \right)+\cdots+t'_k\left(\begin{array}{c}
             x'_k \\
             (x'_k)^{p} \\
             \vdots\\
             (x'_k)^{p^{m-1}}
           \end{array}
  \right),
\end{equation}
where  $x'_1,\ldots,x'_k$ are $\mathbb{F}_p$-linearly independent and $k<m$.
Because the determinant of the coefficient matrix of the system from the first $k$ rows is not zero by Lemma \ref{lem-4.1}, we must have  $t'_i=0$ for all $i$'s, which contradicts with $t'_1(x'_1)^{p^{m-1}}+\ldots+t'_k(x'_{k})^{p^{m-1}}=1$.

If $x_1,\ldots,x_{s}$ are $\mathbb{F}_p$-linearly independent and $s=m$, then the determinant of the coefficient matrix of the system from the first $m$ rows in Eq. (\ref{f-11}) are not zero by Lemma \ref{lem-4.1}. Again we must have $t_i=0$ for all $i$'s, which also  contradicts with $t_1x^{p^{m-1}}_1+\ldots+t_sx^{p^{m-1}}_s=1$. The proof is now complete.
\end{proof}

\section{The girth of linearized Wenger graphs}

In graph theory, the {\it girth} of a graph is the length of a shortest cycle contained in the graph. In \cite{SHS}, Shao et al proved the Wenger graphs  have girth 8, and moreover, if $m \geq 3$, then for any integer $l$ with $l\neq 5, 4 \leq l \leq 2 p$ (where $p$
is the character of the finite field $\mathbb{F}_q$) and any vertex $v$ in the Wenger graph $W_m(q)$, there is a cycle of length $2l$ in $W_m(q)$ passing through the vertex $v$. The existence of the cycles of certain even length plays
an important role in the study of the accurate order of the Tur\'{a}n  number in extremal graph theory. See \cite{bon,chiu,LPS,Murty}.
In this section, we consider the girth of linearized Wenger graphs $L_m(q) = (V, E)$.

Let $P=(p_1,\cdots,p_{m+1}), P'=(p'_1,\cdots,p'_{m+1})$ be two distinct points in $V$. Suppose that $P$ and $P'$ share a common neighbor $L=[l_1,\cdots,l_{m+1}]$, then 
\begin{equation}\label{f-51}
  P-P'=(p_1-p'_1,l_1(p_1-p'_1), l_1(p_1-p'_1)^p,\cdots,l_1(p_1-p'_1)^{p^{m-1}}).
\end{equation}
In other words, $P-P'$ has the form $(u,lu,lu^p,\cdots,lu^{p^{m-1}})$. Conversely, if $P-P'$ has the form $(u,lu,lu^p,\cdots,lu^{p^{m-1}})$ with $u\neq 0$, we show that there exists a unique $L\in V$ such that $L$ is a common neighbor of $P$ and $P'$. Indeed, let $l_1=l$. Since $l_1p_1^{p^{k-2}}-p_k=l_1(p'_1)^{p^{k-2}}-p'_k, k=2,\cdots,m+1$, we can define $l_k=l_1p_1^{p^{k-2}}-p_k,k=2,\cdots,m+1$ and then the point $L=[l_1,\cdots,l_{m+1}]$ is a common neighbor of $P,P'$. Moreover, if both $L=[l_1,\cdots,l_{m+1}]$ and $L'=[l'_1,\cdots,l'_{m+1}]$ are common neighbors of $P,P'$, then by definition, $l_1=l'_1=l$ and $l_k=l'_k=l_1p_1^{p^{k-2}}-p_k=l_1{p'_1}^{p^{k-2}}-p'_k, k=2,\cdots,m+1$. Thus $L=L'$.

We summarize the above discussion as  follows:

{\lem \label{lem-5.1} In the linearized Wenger graph $L_m(q)$, two distinct points $P=(p_1,\cdots,p_{m+1})$ and  $P'=(p'_1,\cdots,p'_{m+1})$ have a common neighbor if and only if $P-P'$ has the form $(u,lu,lu^p,\cdots$, $lu^{p^{m-1}})$ with $u\in \mathbb{F}^*_q,l\in \mathbb{F}_q$. Moreover, if $P-P'$ has the form $(u,lu,lu^p,\cdots,lu^{p^{m-1}})$ with $u\in \mathbb{F}^*_q, l\in \mathbb{F}_q$, then $P,P'$ have a unique common neighbor.}

As a consequence, we have

{\cor  \label{cor-5.2}  There is no cycle of length $4$ in the linearized Wenger graph $L_m(q)$.}
\begin{proof}
If $P_1L_1P_2L_2P_1$ or $L_1P_1L_2P_2L_1$ is a cycle of length $4$ in the linearized Wenger graph, then $L_1,L_2$ are common neighbors of $P_1,P_2$, which is contrary to Lemma \ref{lem-5.1}.\end{proof}

Since the girth of the linearized Wenger graphs is even,  the girth of the linearized Wenger graphs is at least $6$  by Corollary \ref{cor-5.2}.  Furthermore, if $P_1L_1P_2L_2P_3$ $\ldots$ $L_tP_1$ is a cycle of length $2t$ in the linearized Wenger graph $L_m(q)$, then there are elements $u_1,u_2,\ldots,u_t\in \mathbb{F}_q^*$, and $c_1,c_2, \ldots,c_t\in \mathbb{F}_q$ such that
\begin{equation}\label{f-5.2}
\left\{\begin{array}{l}
  P_1-P_2=(u_1,c_1u_1,c_1u_1^p,\cdots,c_1u_1^{p^{m-1}})\\
P_2-P_3=(u_2,c_2u_2,c_2u_2^p,\cdots,c_2u_2^{p^{m-1}})\\
\vdots\\
P_t-P_1=(u_t,c_tu_t,c_tu_t^p,\cdots,c_tu_t^{p^{m-1}})
\end{array}\right .
\end{equation}
and thus
\begin{equation}\label{f-5.3}
  \left\{\begin{array}{l}
           u_1+u_2+\ldots+u_t=0 \\
          c_1u_1+c_2u_2+\ldots+c_tu_t=0 \\
           \vdots\\
               c_1u_1^{p^{m-1}}+c_2u_2^{p^{m-1}}+\ldots+c_tu_t^{p^{m-1}}=0.
         \end{array}
  \right .
\end{equation}

The converse of this result does not hold since $P_1L_1P_2L_2P_3\cdots L_tP_1$  may  not be a cycle. For example,
in linearized Wenger graph $L_1(11)$, choose $P_1=(0,0)$, $P_2=(-1,-1)$, $P_3=(-2,0)$, $P_4=P_1=(0,0)$, $P_5=(-1,-2)$, $P_6=(-2,-8)$,
$L_1=(1,0)$, $L_2=(-1,2)$, $L_3=(0,0)$, $L_4=(2,0)$, $L_5=(6,-4)$, and $L_6=(4,0)$. Then there are
$u_1=u_2=u_4=u_5=1$, $u_3=u_6=-2$, $c_1=1$, $c_2=-1$, $c_3=0$, $c_4=2$, $c_5=6$, $c_6=4$ such that Eq. (\ref{f-5.2}) and (\ref{f-5.3}) hold. However, $P_1L_1\ldots P_6P_1$ is not a cycle in $W_1(11)$.

Therefore,  in order to study cycles of length $2t$ in linearized Wenger graphs, we first  try to solve Eq. (\ref{f-5.2}) and (\ref{f-5.3}). If there are no  $u_i$'s  and $c_i$'s satisfying Eq. (\ref{f-5.2}) and (\ref{f-5.3}), then there is no cycle with length $2t$ in $L_m(q)$. Otherwise, construct $P_1,\ldots,P_t$ and $L_1,\ldots, L_t$ as follows:

 Let $P_i=(p_1^{(i)}, \cdots,p_{m+1}^{(i)}), L_i=[l_1^{(i)}, \cdots,l_{m+1}^{(i)}], i=1,\cdots,t$, where
\begin{eqnarray*}
 p^{(i)}_1-p^{(i+1)}_1=u_i, i=1,2,\ldots,t-1 , p^{(t)}_1-p^{(1)}_1=u_t \\
 l_1^{(i)}=c_i, l_k^{(i)}=l_1^{(i)}(p_1^{(i)})^{p^{k-2}}-p_k^{(i)},k=2,\cdots,m+1.
\end{eqnarray*}

If  both $P_1,\ldots, P_t$ are distinct and $L_1,\ldots,L_t$ are also distinct, then
$P_1L_1P_2L_2P_3$ $\ldots$ $L_tP_1$ is a cycle of length $2t$ in $W_m(q)$. Otherwise, we choose new solutions $u_i$'s and $c_i$'s, and test these new vertices.   If  there are always two $P_i$'s (or two $L_i$'s) which are the same in the above construction for all $u_i$'s and $c_i$'s satisfying Eq. (\ref{f-5.2}) and (\ref{f-5.3}), then  there is no cycle with length $2t$ in $L_m(q)$.


Using the above technique, in the following we give the girth of linearized Wenger graphs.

{\thm \label{thm-5.1}   Let $q=p^e$ and $ m\geq 1$, $e\geq 1$ and $p$ be an odd prime, or  $ m=1$, $e\geq 2$ and $p=2$. Then the girth of the linearized Wenger graph $L_m(q)$  is 6.}

\begin{proof} Case 1. $ m\geq 1$, $e\geq 1$ and $p$ is an odd prime. By Corollary \ref{cor-5.2}, it is enough to construct a cycle with length $6$ in this case. Indeed, let $u_1=u_2=1,u_3=-2$, $c_1=1$, $c_2=-1$, $c_3=0$, $P_1=(0,0,\ldots,0)$, $P_2=(-1,-1,\ldots, -1)$, $P_3=(-2,0,\ldots,0)$, $L_1=[1,0,\ldots,0]$, $L_2=[-1,2,2,\ldots,2]$, $L_3=[0,0,\ldots,0]$. Then $P_1L_1P_2$ $L_2$ $P_3L_3P_1$ is a cycle with length $6$.

Case 2.  $e\geq 2$, $m=1$ and $p=2$.  For an element $\beta\in \mathbb{F}^{*}_q$ and ${\rm tr}(\beta)=0$, there exists some $\alpha\in \mathbb{F}^{*}_q$ such that $\alpha^2+\alpha=\beta$.
Put $u_1=\alpha^2$, $u_2=\alpha$, $u_3=\beta$,  $c_1=0$, $c_2=\alpha^{-1}\beta $  and $c_3=1$. One can construct a cycle $P_1L_1P_2L_2P_3L_3P_1$ of length $6$, where $P_1=(0,0)$,
$P_2=(\alpha^2,0)$, $P_3=(\beta,\beta)$, $L_1=[0,0]$, $L_2=[\alpha^{-1}\beta,\alpha\beta]$ and $L_3=[1,0]$.
\end{proof}

{\thm \label{thm-5.2}  Let $q=p^e$,  $p=2$ and either $e=m=1$ or $e\geq 1$, $m\geq 2$. Then the girth of
the linearized Wenger graph $L_m(q)$ is $8$. }

\begin{proof}
First we need to show that there is no  cycle of length $6$ in $L_{m}(q)$ in these two cases. For the case of $e=1$ and $p=2$, there is no $u_i\in \mathbb{F}^{*}_q$, $1\leq i\leq 3$, such that Eq (\ref{f-5.3}) holds. Hence there is no cycle with length $6$ in this case.  Assume that there is a cycle $P_1L_1P_2L_2P_3L_3P_1$ of length $6$ in $L_m(q)$ for the case of $e\geq 2$, $m\geq 2$ and $p=2$. Then there are elements $u_1,u_2,u_3\in \mathbb{F}^{*}_q$,  $c_1,c_2,c_3\in \mathbb{F}_q$ such that Eq (\ref{f-5.2}) and (\ref{f-5.3}) hold.

Eliminating $c_1$ among two successive equations of the last $m-1$ equations in  Eq. (\ref{f-5.3}), we get
\begin{equation}\label{f-5.4}
  \left\{\begin{array}{l}
           u_1+u_2+u_3=0 \\
          c_1u_1+c_2u_2+c_3u_3=0 \\
        c_2(u^{2}_2-u_2u_1)+c_3(u^{2}_3-u_3u_1)=0 \\
           \vdots\\
               c_2(u_2^{2^{m-1}}-u_{2}^{2^{m-2}}u_{1}^{2^{m-2}})+c_3(u_3^{2^{m-1}}-u_{3}^{2^{m-2}}u_{1}^{2^{m-2}})=0.
         \end{array}
  \right .
\end{equation}
Further simplifying Eq. (\ref{f-5.4}) by using $u_1+u_2+u_3=0$ and $u_1, u_2, u_3 \in \mathbb{F}_q^*$, we get
\begin{equation}\label{f-5.5}
  \left\{\begin{array}{l}
           u_1+u_2+u_3=0 \\
          c_1u_1+c_2u_2+c_3u_3=0 \\
        c_2+c_3=0 \\
           \vdots\\
               c_2+c_3=0.
         \end{array}
  \right .
\end{equation}
 Therefore, by symmetry,  Eq. (\ref{f-5.3})  has only the solution $c_1=c_2=c_3$. Then we have $L_1=L_3$ since they share the common vertex $P_1$,  which contradicts to the earlier assumption.

In the following we can construct a cycle $P_1L_1P_2L_2\ldots L_4P_1$ in both cases:
Put $u_1=u_2=u_3=u_4=1$ and $c_1=c_3=0$, $c_2=c_4=1$. Let $P_1=(0,0,0,\ldots,0)$,  $P_2=(1,0,0,\ldots,0)$, $P_3=(0,1,1,\ldots,1)$,
$P_4=(1,1,1, \ldots,1)$, $L_1=[0,0,0,\ldots,0]$, $L_2=[1,1,1, \ldots,1]$, $L_3=[0,1,1,\ldots,1]$, $L_4=[1,0,0,\ldots, 0]$. Then it is straightforward to check $P_1L_1P_2L_2\ldots L_4P_1$ is indeed a cycle of length $8$. Hence we complete the proof.
\end{proof}

\section{Open Problems}
There are several open problems about linearized Wenger graphs.
First finding an explicit formula for the eigenvalue multiplicities $n_{p^i}$'s of the linearized Wenger graphs when $m<e$ is an open problem. Constructing even cycles with specific length in linearized Wenger graphs is also interesting. In addition, it would be desirable to find new classes of  $f_k(x)$ such that the explicit spectrum of these new types of Wenger graphs can be determined by Theorem 2.2.



\begin{thebibliography}{amsplain}






\bibitem{babai} L. Babai, Spectra of Cayley graphs, J. Combin. Theory Ser. B, {\bf 27} (1979) 180-189.
\bibitem{boll}B. Bollobas, Modern Graph Theory, Springer-Verlag New York, Inc., 1998.
\bibitem{bon} J. Bondy, M. Simonovits, Cycles of even length in graphs. J. Combin. Theory (Series B), {\bf 16}
(1974) 97-105.
\bibitem{chiu}P. Chiu, Cubic Ramanujan graphs, Combinatorica, {\bf 12} (1992) 275-285.
\bibitem{cll}S. M. Cioab$\rm{\breve{a}}$, F. Lazebnik, W. Li, On the spectrum of Wenger graphs, J. Combin. Theory, Ser. B, {\bf 107} (2014) 132-139.

\bibitem{DLW} V. Dmytrenko, F.  Lazebnik and J. Williford, On monomial graphs of girth eight,
Finite Fields Appl., {\bf 13}(4) (2007) 828-842.



\bibitem{HLW} S. Hoory, N. Linial, A. Wigderson, Expanders and their applications, Bull. Amer. Math. Soc. {\bf 43} (2006) 439-561.
\bibitem{jukna}S. Jukna, Extremal Combinatorics, Texts in Theoretical Computer Science. Springer-Verlag Berlin Heidelberg, 2011.
\bibitem{lidl}R. Lidl, H. Niederriter,  Finite Fields, Encyclopedia Math.
Appl. Vol. {\bf 20}, Addison-Wesley, Reading, 1983.
\bibitem{ling}S. Ling, L.J. Qu,  A note on linearized polynomials and the dimension of their
kernels, Finite Fields Appl., 18 (2012) 56-62.
\bibitem{lu}F. Lazebnik, V.Ustimenko, New examples of graphs without small cycles and of large size, European J. Combin., {\bf 14} (1993) 445-460.
\bibitem{lu2}F. Lazebnik, V. Ustimenko, Explicit construction of graphs with arbitrary large girth and of large size, Discrete Appl. Math., {\bf 60} (1997) 275-284.
\bibitem{LV}  F. Lazebni, R. Viglione,  An infinite series of regular edge- but not vertex transitive graphs, J. Graph Theory {\bf 41} (2002) 249-258.
    \bibitem{lovasz} L. Lov$\acute{a}$sz, Spectra of graphs with transitive groups, Period. Math. Hungar. 6 (1975)  191-195
\bibitem{LPS} A. Lubotzky, R. Phillips, P. Sarnak, Ramanujan graphs, Combinatorica, {\bf 8}(3) (1988) 261-277.
\bibitem{MM}K. Mellinger, D. Mubayi, Constructions of bipartite graphs from finite geometries, J. Graph Theory, {\bf 49}(1) (2005) 1-10.
\bibitem{Murty}M.R. Murty, Ramanujan graphs, J. Ramanujan Math. Soc., {\bf 23} (2003) 33-52.

\bibitem{SHS} J.-Y. Shao, C.-X. He, H.-Y.Shan, The existence of even cycles with specific lengths in Wenger's graph, Acta Math. Appl. Sin. Engl. Ser., {\bf 24} (2008) 281-288.

\bibitem{Viglione0} R. Viglione, Properties of some algebraically defined graphs, PhD thesis, University of Delaware, 2002.

\bibitem{Viglione} R. Viglione, On the diameter of Wenger graphs, Acta Appl. Math. {\bf 104} (2008) 173-176.
\bibitem{Wenger} R. Wenger, Extremal graphs with no $C^4$'s, $C^6$'s, or $C^{10}$'s, J. Combin. Theory Ser. B, {\bf 52}(1) (1991) 113-116.




























\end{thebibliography}
\end{document}